\documentclass[10pt]{article}
\usepackage{epsfig,amssymb,amsmath}
\usepackage{palatino}
\usepackage{verbatim}
\usepackage{enumerate}

\def\baa {\begin{eqnarray*}}
\def\eaa {\end{eqnarray*}}

\def \la {\lambda}
\def \al {\alpha}
\def \be {\beta}

\def\la{\lambda}

\newtheorem{lemma}{Lemma}[section]
\newtheorem{proposition}[lemma]{Proposition}
\newtheorem{corollary}[lemma]{Corollary}
\newtheorem{theorem}[lemma]{Theorem}

\newtheorem{problem}[lemma]{Problem}

\def\bc  {\begin{comment}}
\def\ec  {\end{comment}}
\makeatletter
\@addtoreset{equation}{section}
\makeatother
\def\proof{\medskip\noindent{\bf Proof.} }
\def\qed{\hfill $\Box$}

\newcommand {\ds} {\displaystyle}

\begin{document}

\title{Some inequalities for Chebyshev polynomials}

\author{Geno Nikolov}

\date{}
\maketitle

%******************************************************************
\begin{abstract}
Askey and Gasper (1976) proved a trigonometric inequality which
improves another trigonometric inequality found by M. S. Robertson
(1945). Here these inequalities are reformulated in terms of the
Chebyshev polynomial of the first kind $T_n$ and then put into a
one-parametric family of inequalities. The extreme value of the
parameter is found for which these inequalities hold true. As a step
towards the proof of this result we establish the following
complement to the finite increment theorem specialized to
$T_n^{\prime}$:
$$
T_n^{\prime}(1)-T_n^{\prime}(x)\geq
(1-x)\,T_n^{\prime\prime}(x)\,,\qquad x\in [0,1]\,.
$$
By a known expansion formula, this property is extended for the
class of ultraspherical polynomials $P_n^{(\la)}$, $\la\geq 1$.
\end{abstract}

\textbf{MSC 2010:} 41A17
\smallskip

\textbf{Key words and phrases:} Positive trigonometric sums,
Chebyshev polynomials
%******************************************************************
\section{Introduction and statement of the results}
Positive trigonometric sums play an important role in Harmonic
Analysis, Orthogonal Polynomials, Approximation Theory and many
other branches of mathematics. In this note we discuss two
trigonometric inequalities which appear in the book of R. Askey
\cite{RA1975}. The first one is inequality (1.29) in \cite{RA1975},
which reads as
\begin{equation}\label{e1.1}
\frac{\sin (n-1)\theta}{(n-1)\sin \theta}-\frac{\sin
(n+1)\theta}{(n+1)\sin \theta}\leq\frac{4n}{n^2-1}\Big[1-\frac{\sin
n\theta}{n\sin \theta}\Big]\,,\qquad 0\leq\theta\leq\pi
\end{equation}
(actually, (1.29) appears in \cite{RA1975} as a strict inequality
and under the assumption $0<\theta<\pi$). It was proved by M. S.
Robertson \cite{MR1945} while studying the coefficients of univalent
functions. The second one was proved by Askey and Gasper
\cite{AG1976} and reads as
\begin{equation}\label{e1.2}
\frac{\sin (n-1)\theta}{(n-1)\sin \theta}-\frac{\sin
(n+1)\theta}{(n+1)\sin \theta}\leq
\frac{(3+\cos\theta)n}{n^2-1}\Big[1-\frac{\sin n\theta}{n\sin
\theta}\Big]\,,\qquad 0\leq\theta\leq\pi\,.
\end{equation}
This is inequality (8.17) in \cite{RA1975}, and as Askey wrote, it
is sharper than \eqref{e1.1}. For more information on positive
trigonometric sums and positive finite linear combinations of
classical orthogonal polynomials we refer to \cite{RA1975, DM1998,
DM2001} and the references therein.

We find it more convenient to reformulate inequalities \eqref{e1.1}
and \eqref{e1.2} in terms of the Chebyshev polynomials of the first
kind. Let us recall that the $m$-th Chebyshev polynomial of the
first kind $T_m$, $m\in \mathbb{N}_0$, is defined by
$$
T_m(x)=\cos m\theta,\qquad x=\cos\theta\in [-1,1],\ \ \theta\in
[0,\pi],
$$
and its derivative is
$$
T_m^{\prime}(x)=m\,\frac{\sin m\theta}{\sin \theta},\qquad
x=\cos\theta\,.
$$
Using
\begin{eqnarray*}
&&\frac{\sin (n-1)\theta}{\sin\theta}=\frac{\sin
n\theta\cos\theta-\cos
n\theta\sin\theta}{\sin\theta}=\frac{xT_n^{\prime}(x)}{n}-T_n(x)\,,\\
&&\frac{\sin (n+1)\theta}{\sin\theta}=\frac{\sin
n\theta\cos\theta+\cos
n\theta\sin\theta}{\sin\theta}=\frac{xT_n^{\prime}(x)}{n}+T_n(x)\,,
\end{eqnarray*}
we find that inequality  \eqref{e1.1} of Robertson is equivalent to
the inequality
\begin{equation}\label{e1.3}
f_1(x):=T_n(x)+2-\frac{x+2}{n^2}\,T_n^{\prime}(x)\geq 0\,,\qquad
n\geq 2,\ \ x\in [-1,1]\,,
\end{equation}
while the Askey-Gasper inequality \eqref{e1.2} is equivalent to the
inequality
\begin{equation}\label{e1.4}
f_2(x):=T_n(x)+\frac{x+3}{2}-\frac{3(x+1)}{2n^2}\,T_n^{\prime}(x)
\geq 0\,,\qquad n\geq 2,\ \  x\in [-1,1]\,.
\end{equation}
Since $\max_{x\in [-1,1]}T_n^{\prime}(x)=T_n^{\prime}(1)=n^2$, we
have
$$
f_1(x)-f_2(x)=\frac{1-x}{2n^2}\,\big(n^2-T_n^{\prime}(x)\big) \geq
0,\qquad x\in [-1,1]\,,
$$
hence inequality \eqref{e1.3} is a consequence of inequality
\eqref{e1.4}. Thus we naturally arrive at the following
\begin{problem}\label{p1}
Find the largest constant $a=a(n)\geq 0$, $\,n\geq 2$, such that
$$
g_n(a;x):=T_n(x)+2-\frac{x+2}{n^2}\,T_n^{\prime}(x)-a\,\frac{1-x}{n^2}\,
\big(n^2-T_n^{\prime}(x)\big)\geq 0,\qquad x\in [-1,1]\,.
$$
\end{problem}
The cases $n=2,\;3$ are trivial: from $\,T_2(x)=2x^2-1\,$ and
$\,T_3(x)=4x^3-3x\,$ one finds
$$
g_2(a;x)=(1-a)(1-x)^2\,,\qquad
g_3(a;x)=\frac{4}{3}(2-a)(1-x)^2(1+x)\,,
$$
whence $a(2)=1$ and $a(3)=2$.

We assume henceforth that $n\geq 4$. By the Askey--Gasper inequality
\eqref{e1.2}, $g_n(1/2,x)=f_2(x)\geq 0$, $x\in [-1,1]$, and
therefore $a(n)\geq 1/2$. We show below that $a(n)$ cannot be
essentially larger than $1/2$. Indeed, let
$$
x_k=\cos \frac{k\pi}{n}, \qquad k=1,\ldots,n-1,
$$
be the zeros of $T_n^{\prime}$. Then $\,T_n(x_k)=(-1)^{k}\,$ and
$$
g_n(a;x_k)=2+(-1)^k-a(1-x_k),\qquad k=1,\ldots,n-1\,.
$$
The condition that $\,g_n(a;x_{n-1})\geq 0$ in the case of even $n$
and $\,g_n(a;x_{n-2})\geq 0$ in the case of odd $n$ implies
respectively
$$
a\leq\frac{1}{1-x_{n-1}}=\frac{1}{1+x_{1}}
=\frac{1}{1+\cos\frac{\pi}{n}},\qquad n - \text{ even,}
$$
and
$$
a\leq\frac{1}{1-x_{n-2}}=\frac{1}{1+x_{2}}
=\frac{1}{1+\cos\frac{2\pi}{n}},\qquad n - \text{ odd.}
$$
Since both upper bounds for $a$ tend to $1/2$ as $n$ grows, it
follows that $a=1/2$ is the best possible (the largest) absolute
constant, ensuring that $g_n(a;x)\geq 0$ for every $x\in [-1,1]$ and
all $n\in \mathbb{N}$, $n\geq 2$. In this sense, the Askey-Gasper
inequality \eqref{e1.4} is the best possible.

It turns out that the upper bounds for $a(n)$ found above actually
provide the solution to Problem~\ref{p1}. Specifically, we prove the
following statement.
\begin{theorem}\label{t1}
Let $n\in \mathbb{N}$, $n\geq 4$. Then
\begin{equation}\label{e1.5}
T_n(x)+2-\frac{x+2}{n^2}\,T_n^{\prime}(x)-a(n)\,\frac{1-x}{n^2}\,
\big(n^2-T_n^{\prime}(x)\big)\geq 0\,,\quad  x\in [-1,1]\,,
\end{equation}
where
$$
a(n)=\begin{cases}\ds{\frac{1}{1+\cos\frac{\pi}{n}}},& \text{ if }\
n\ \text{ is even},\vspace*{1ex}\\
\ds{\frac{1}{1+\cos\frac{2\pi}{n}}},&\text{ if }\ n \ \text{ is
odd}.
\end{cases}
$$
The constant $a(n)$ is the best possible in the sense that
\eqref{e1.5} fails for any larger constant. The equality in
\eqref{e1.5} is attained only at $x=1$ and $x=-\cos\frac{\pi}{n}$ if
$n$ is even, and at $x=\pm 1$ and $x=-\cos\frac{2\pi}{n}$ if $n$ is
odd.
\end{theorem}

\begin{figure}[htp]
\centering
\includegraphics[scale=0.65,clip]{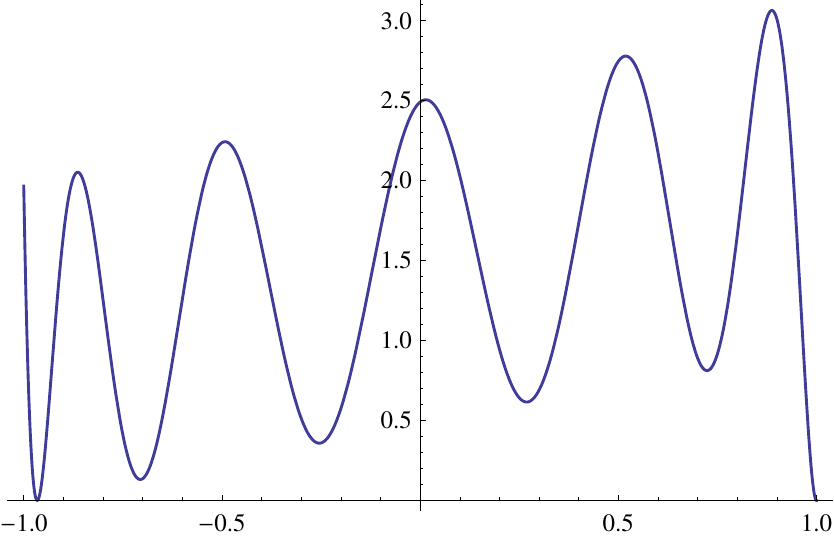} \hspace*{0.5cm}
\includegraphics[scale=0.65,clip]{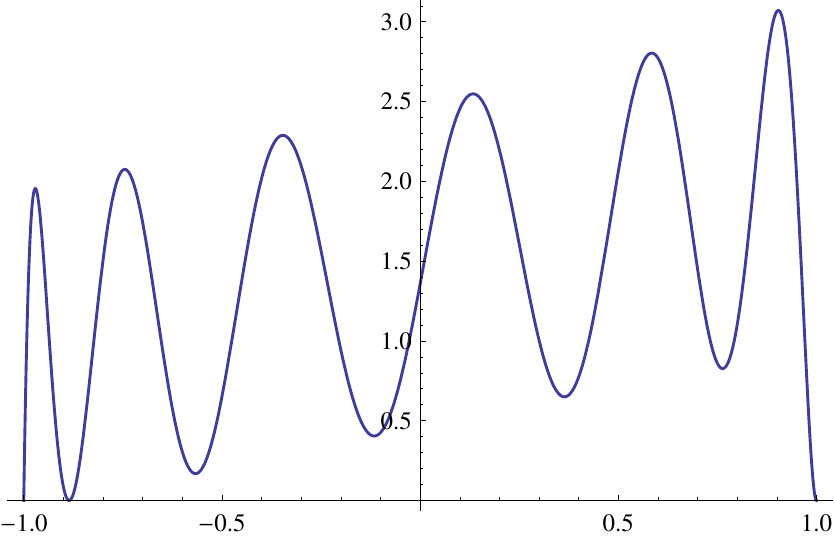}
{\small \caption{The graphs of
$T_n(x)+2-\frac{x+2}{n^2}\,T_n^{\prime}(x)-a(n)\,\frac{1-x}{n^2}\,
\big(n^2-T_n^{\prime}(x)\big)$ for $n=12$ (left) and $n=13$
(right).}} \label{fig1}
\end{figure}

The typical behavior of the function in \eqref{e1.5} in the cases of
even and odd $n$ is shown in Figure~1. The graphs suggest that  in
the interval  $[0,1]$  this function could be non-negative for a
larger constant $a$ than the one specified in Theorem~\ref{t1}. We
show that this is indeed the case by proving that the
non-negativeness in $[0,1]$ persists with $a=1$.

\begin{theorem}\label{t2}
Let $n\in \mathbb{N}$, $n\geq 3$. Then
\begin{equation}\label{e1.6}
T_n(x)+x+1-\frac{2x+1}{n^2}\,T_n^{\prime}(x)\geq 0\,,\quad x\in
[0,1]\,,
\end{equation}
or, equivalently,
\begin{equation}\label{e1.7}
T_n^{\prime}(1)-T_n^{\prime}(x)\geq
(1-x)\,T_n^{\prime\prime}(x)\,,\qquad x\in [0,1]\,.
\end{equation}
The equality in \eqref{e1.6}-\eqref{e1.7} occurs only for $x=1$ and
if $n\equiv 2$ $(\!\!\!\mod 4)$, for $x=0$.
\end{theorem}
For $n=0,\,1,\,2,$ \eqref{e1.7} becomes an identity. Inequality
\eqref{e1.7} provides an interesting complement to the finite
increment formula $T_n^{\prime}(1)-T_n^{\prime}(x)=
(1-x)\,T_n^{\prime\prime}(\xi)$, $\xi\in (x,1)$. Moreover,
\eqref{e1.7} implies a similar property of the ultraspherical
polynomials $P_n^{(\la)}$, $\la\geq 1$.

\begin{corollary}\label{c1}
For every ultraspherical polynomial $P_n^{(\la)}$, $\la\geq 1$,
there holds
$$
P_n^{(\la)}(1)-P_n^{(\la)}(x)\geq
\frac{d}{dx}\,\big\{P_n^{(\la)}(x)\big\}\,(1-x)\,,\qquad x\in
[0,1]\,.
$$
\end{corollary}

Corollary~\ref{c1} easily follows from
$T_n^{\prime}=n\,P_{n-1}^{(1)}$ and the fact that if $\mu\geq \la$,
then $P_n^{(\mu)}$ is represented as a linear combination of
$\{P_m^{(\la)}\}_{m=0}^{n}$ with non-negative coefficients. There
are examples showing that Corollary~\ref{c1} is not true if $\la<1$,
the case $\la=0$ is particularly easy to verify. A challenging
problem is to characterize all pairs of parameters $(\al,\be)$
ensuring similar inequality for the Jacobi polynomials
$P_n^{(\al,\be)}$.
\medskip

The paper is organized as follows. In the next section we propose a
short elementary proof of the Askey-Gasper inequality \eqref{e1.4}
(and thereby of the Robertson inequality \eqref{e1.3}). In Section~3
we present a proof of Theorem~\ref{t2}. The proof of
Theorem~\ref{t1} is given in Section~4.

%==================================================================
\section{Proof of the Askey--Gasper inequality (1.4)}
%==================================================================
\setcounter{equation}{0} We prove the following statement:
\begin{proposition}\label{pp1}
Let $n\in \mathbb{N}$,  $n\geq 2$. Then
\begin{equation}\label{e2.1}
f_2(x)=T_n(x)+\frac{x+3}{2}-\frac{3(x+1)}{2n^2}\,T_n^{\prime}(x)\geq
0,\qquad x\in [-1,1].
\end{equation}
The equality in \eqref{e2.1} is attained only at $x=1$ and if $n$ is
odd, at $x=-1$.
\end{proposition}

\proof From $\,T_2(x)=2x^2-1\,$ and $\,T_3(x)=4x^3-3x\,$ we find
\begin{eqnarray*}
&& f_2(x)=\frac{1}{2}\,(1-x)^2,\ \ n=2,\\
&& f_2(x)=2(1-x)^2(1+x),\ \ n=3,
\end{eqnarray*}
hence Proposition~\ref{pp1} is true for $\,n=2,\,3$, and we assume
henceforth $n\geq 4$.

We shall use in this and in the next section the differential
equation satisfied by $T_n$,
\begin{equation}\label{e2.2}
(1-x^2)\,T_n^{\prime\prime}(x)-x\,T_n^{\prime}(x)+n^2\,T_n(x)=0\,,
\end{equation}
as well as the identity
\begin{equation}\label{e2.3}
n^2\,\big[T_n(x)\big]^2+(1-x^2)\,\big[T_n^{\prime}(x)\big]^2=n^2\,.
\end{equation}
Denote by $\tau$ the largest zero of $T_n^{\prime\prime}$. \medskip

\emph{Case 1: $x\in (\tau,1]$}. Using \eqref{e2.2}, we rewrite $f_2$
in the form
$$
f_2(x)=\frac{x+3}{2n^2}\,\big(T_n^{\prime}(1)-T_n^{\prime}(x)\big)
-\frac{1}{n^2}\,(1-x^2)T_n^{\prime\prime}(x)\,.
$$
We observe that in this case the inequality $f_2(x)\geq 0$ is
equivalent to
\begin{equation}\label{e2.4}
\frac{x+3}{2(x+1)}\,\frac{1}{1-x}\,\int_{x}^{1}T_n^{\prime\prime}(u)\,du
\geq T_n^{\prime\prime}(x)\,,\qquad x\in (\tau,1]\,.
\end{equation}
Since $T_n^{\prime\prime}$ is positive and monotonically increasing
in $(\tau,1]$,
$$
\frac{1}{1-x}\,\int_{x}^{1}T_n^{\prime\prime}(u)\,du\geq
T_n^{\prime\prime}(x)
$$
and \eqref{e2.4} is a consequence of the inequality
$$
\frac{x+3}{2(x+1)}\geq 1\,,
$$
which is obviously true for $x\in (\tau,1]$. Notice in the last two
inequalities the equality holds only when $x=1$. Hence, $f_2(x)\geq
0$ for $x\in (\tau,1]$, and the equality is attained only for $x=1$.
\medskip

\emph{Case 2: $x\in [-1,\tau]$}. In this case we rewrite inequality
$f_2(x)\geq 0$ in the form
\begin{equation}\label{e2.5}
1+T_n(x)+\frac{1+x}{2}\geq\frac{3(1+x)}{2n^2}\,T_n^{\prime}(x),
\qquad x\in [-1,\tau]\,.
\end{equation}
The left-hand side of \eqref{e2.5} is non-negative for $x\in
[-1,\tau]$. On the other hand, the right-hand side of \eqref{e2.5}
is non-positive for $x\in [\cos\frac{2\pi}{n},\cos\frac{\pi}{n}]$ as
$T_n^{\prime}(x)\leq 0$ therein. Since $\tau\in
(\cos\frac{2\pi}{n},\cos\frac{\pi}{n})$, \eqref{e2.5} will be proved
if we show that
$$
1+T_n(x)+\frac{1+x}{2}\geq\frac{3(1+x)}{2n^2}\,T_n^{\prime}(x),
\qquad x\in \big[-1,\cos\frac{2\pi}{n}\big]\,,
$$
and it suffices to prove the inequality
\begin{equation}\label{e2.6}
\Big(1+T_n(x)+\frac{1+x}{2}\Big)^2
\geq\frac{9(1+x)^2}{4n^4}\,\big[T_n^{\prime}(x)\big]^2, \qquad x\in
\big[-1,\cos\frac{2\pi}{n}\big]\,.
\end{equation}
We estimate the left-hand side of \eqref{e2.6} by the arithmetic
mean - geometric mean inequality:
\begin{equation}\label{e2.7}
\Big(1+T_n(x)+\frac{1+x}{2}\Big)^2\geq 2(1+x)(1+T_n(x))\,,
\end{equation}
and apply identity \eqref{e2.3} to express
$\big[T_n^{\prime}(x)\big]^2$ in the right-hand side of
\eqref{e2.6},
$$
\big[T_n^{\prime}(x)\big]^2=\frac{n^2(1-T_n(x))(1+T_n(x))}{1-x^2}\,.
$$
It follows from \eqref{e2.7} that \eqref{e2.6} will hold for $x\in
\big[-1,\cos\frac{2\pi}{n}\big]$  if
\begin{equation}\label{e2.8}
2(1+x)(1+T_n(x))\geq\frac{9(1+x)}{4n^2(1-x)}\,(1+T_n(x))(1-T_n(x))\,.
\end{equation}
Since $(1+x)(1+T_n(x))\geq 0$ and $1-T_n(x)\leq 2$, the above
inequality will be certainly true if
$$
n^2(1-x)\geq \frac{9}{4}\,,\qquad x\in
\big[-1,\cos\frac{2\pi}{n}\big]\,.
$$
To see that the above inequality is true, we make use of
$\sin\al>\frac{2}{\pi}\,\al$, $\al\in (0,\pi/2)$. For $x\in
\big[-1,\cos\frac{2\pi}{n}\big]$,
$$
n^2(1-x)\geq
n^2\Big(1-\cos\frac{2\pi}{n}\Big)=2n^2\sin^2\frac{\pi}{n}>
2n^2\,\Big(\frac{2}{\pi}\,\frac{\pi}{n}\Big)^2=8>\frac{9}{4}\,.
$$
Thus, the inequality $f_2(x)\geq 0$ is proved in Case 2, too, and it
remains to check when the equality is attained. Tracing backward our
proof, we see that the equality in \eqref{e2.8} holds only if either
$x=-1$ or $T_n(x)+1=0$, while the equality in \eqref{e2.7} holds
only if $T_n(x)+1=(1+x)/2$. Both conditions imply that $x=-1$ and
$T_n(-1)=-1$, and the latter holds if and only if $n$ is odd. \qed
%==================================================================
\section{Proof of Theorem 2}
%==================================================================
\setcounter{equation}{0}

Let us set
\begin{eqnarray*}
\varphi_n(x)&:=&T_n(x)+x+1-\frac{2x+1}{n^2}\,T_n^{\prime}(x)\,,\\
\psi_n(x)&:=&T_n^{\prime}(1)-T_n^{\prime}(x)-(1-x)\,T_n^{\prime\prime}(x)\,.
\end{eqnarray*}
From $T_n^{\prime}(1)=n^2$ and the differential equation
\eqref{e2.2} we have
\begin{equation*}
\begin{split}
\varphi_n(x)&=\frac{x+1}{n^2}\,\big(n^2-T_n^{\prime}(x)\big)
-\frac{1}{n^2}\,\big(x\,T_n^{\prime}(x)-n^2\,T_n(x)\big)\\
&=\frac{x+1}{n^2}\,\big(T_n^{\prime}(1)-T_n^{\prime}(x)\big)-\frac{1-x^2}{n^2}
\,T_n^{\prime\prime}(x)\\
&=\frac{x+1}{n^2}\,\psi_n(x)\,,
\end{split}
\end{equation*}
therefore inequalities  \eqref{e1.6} and \eqref{e1.7}, i.e.,
$\varphi_n(x)\geq 0$ and $\psi_n(x)\geq 0$, $x\in [0,1]$,  are
equivalent.

From
$$
T_n(x)=\cos n\theta,\ \ T_n^{\prime}(x)=n\,\frac{\sin
n\theta}{\sin\theta}\,,\qquad x=\cos\theta\,,\ \ \theta\in
[0,\pi]\,,
$$
we find
$$
\varphi_n(0)=\begin{cases} 2,& n\equiv 0\ (\!\!\!\!\mod 4)\\
1-1/n,& n\equiv 1\ (\!\!\!\!\mod 4)\\
0,& n\equiv 2\ (\!\!\!\!\mod 4)\\
1+1/n,& n\equiv 3\ (\!\!\!\!\mod 4)\,,
\end{cases}
$$
hence $\varphi_n(0)\geq 0$, with the equality holding only if
$n=4k+2$, and we may assume further $x\in (0,1]$.

If $\tau>0$ is the largest zero of $T_n^{\prime\prime}$, then, by
the same argument as in \emph{Case~1} in the preceding section, we
obtain
$$
\psi_n(x)=(1-x)\Big\{\frac{1}{1-x}\int_{x}^{1}T_n^{\prime\prime}(u)\,du-
T_n^{\prime\prime}(x)\Big\}\geq 0,\qquad x\in [\tau,1]\,,
$$
with the equality holding only for $x=1$, so we may restrict our
consideration to the case $x\in (0,\tau)$. Furthermore,
$T_n^{\prime}(x)\leq 0$ for $x\in
\big[\cos\frac{2\pi}{n},\cos\frac{\pi}{n}\big]$, and then obviously
$\varphi_n(x)>0$ for $x\in
\big[\cos\frac{2\pi}{n},\cos\frac{\pi}{n}\big]$. Since $\tau$ lies
in this interval, it remains to prove either of the inequalities
$\varphi_n(x)>0$ and $\psi_n(x)>0$ when $x\in
(0,\cos\frac{2\pi}{n})$. There is nothing to prove if $n=4$, so we
assume $n\geq 5$.  It suffices to show that $\varphi_n(t)>0$ (or
$\psi_n(t)>0$) for every critical point $t$ of $\psi_n$ (i.e. zero
of $\psi_n^{\prime}$) in the interval $(0,\cos\frac{2\pi}{n})$.
Since
$$
\psi_n^{\prime}(x)=(x-1)T_n^{\prime\prime\prime}(x)\,,
$$
the critical points of $\psi_n$ in $(0,\cos\frac{2\pi}{n})$ are
zeros of $T_n^{\prime\prime\prime}$.

Let $t\in (0,\cos\frac{2\pi}{n})$ be a zero of
$T_n^{\prime\prime\prime}$. From
$$
\cos\frac{2\pi}{n}=1-2\,\sin^2\frac{\pi}{n}
<1-2\,\Big(\frac{2}{\pi}\,\frac{\pi}{n}\Big)^2=1-\frac{8}{n^2},
$$
we conclude that
\begin{equation}\label{e3.1}
0<t<1-\frac{8}{n^2}\,.
\end{equation}
Since $y=T_n^{\prime}(x)$ satisfies the differential equation
$$
(1-x^2)y^{\prime\prime}-3x\,y^{\prime}+(n^2-1)\,y=0
$$
(this can be seen e.g., by differentiating \eqref{e2.2}), and
$y^{\prime\prime}(t)=0$, we have
$$
T_n^{\prime\prime}(t)=\frac{n^2-1}{3t}\,T_n^{\prime}(t)\,.
$$
Now, \eqref{e2.2} with $x=t$ yields
$$
(1-t^2)\,\frac{n^2-1}{3t}\,T_n^{\prime}(t)-t\,T_n^{\prime}(t)+n^2\,T_n(t)=0\,,
$$
whence
\begin{equation}\label{e3.2}
\frac{T_n^{\prime}(t)}{n^2}=-\frac{3t\,T_n(t)}{n^2-1-(n^2+2)t^2}\,.
\end{equation}
Let us point out that \eqref{e3.1} implies $n^2-1-(n^2+2)t^2>0$,
since
$$
t^2<t<1-\frac{8}{n^2}<1-\frac{3}{n^2+2}=\frac{n^2-1}{n^2+2}\,.
$$
Replacing $T_n^{\prime}(t)/n^2$ with the right-hand side of
\eqref{e3.2}  in $\varphi_n(t)$, we obtain
\begin{equation*}
\begin{split}
\varphi_n(t)&=t+1+\Big(1+\frac{3t(2t+1)}{n^2-1-(n^2+2)t^2}\Big)\,T_n(t)\\
&=(t+1)\Big\{1+\frac{n^2-1-(n^2-4)t}{n^2-1-(n^2+2)t^2}\,T_n(t)\Big\}\\
&\geq (t+1)\,\Big\{1-\frac{n^2-1-(n^2-4)t}{n^2-1-(n^2+2)t^2}\Big\}
\end{split}
\end{equation*}
(we have used that the factor in front of $T_n(t)$ in the curly
brackets is positive and $T_n(t)\geq -1$). Hence, to prove
$\varphi_n(t)>0$ it suffices to show that
$$
1-\frac{n^2-1-(n^2-4)t}{n^2-1-(n^2+2)t^2}>0\,,
$$
which is equivalent to
$$
\frac{t\,\big[n^2-4-(n^2+2)t\big]}{n^2-1-(n^2+2)t^2}>0\,.
$$
This inequality is true, since the numerator in the left-hand side
is positive. Indeed, from \eqref{e3.1} we have
$$
t<1-\frac{8}{n^2}<1-\frac{6}{n^2+2}=\frac{n^2-4}{n^2+2}.
$$
The proof of Theorem~\ref{t2} is complete.\qed

%==================================================================
\section{Proof of Theorem 1}
%==================================================================
\setcounter{equation}{0} Recall that
$$
f_1(x)=T_n(x)+2-\frac{x+2}{n^2}\,T_n^{\prime}(x),
$$
and set
$$
f_3(x):=\frac{1}{n^2}\,(1-x)(n^2-T_n^{\prime}(x)),
$$
$$
F_a(x):=(1+a)f_1(x)-f_3(x)\,.
$$
Clearly, Theorem~\ref{t1} is equivalent to the following statement:
\begin{theorem}\label{t3}
Let $n\in \mathbb{N}$, $n\geq 4$ and
\begin{equation}\label{e4.1}
a=\begin{cases}\cos\frac{\pi}{n},& \text{ if }\
n\ \text{ is even},\vspace*{1ex}\\
\cos\frac{2\pi}{n},&\text{ if }\ n \ \text{ is odd}.
\end{cases}
\end{equation}
Then
\begin{equation}\label{e4.2}
F_a(x)\geq 0\,,\qquad x\in [-1,1]\,.
\end{equation}
The equality in \eqref{e4.2} occurs only for $x=1$ and
$x=-\cos\frac{\pi}{n}$ if $n$ is even, and for $x=\pm 1$ and
$x=-\cos\frac{2\pi}{n}$ if $n$ is odd. For every constant $a$,
smaller than the one specified in \eqref{e4.1}, inequality
\eqref{e4.2} fails to hold.
\end{theorem}

\proof According to Robertson's inequality, $f_1(x)\geq 0$ in
$[-1,1]$, therefore if $a_1>a_2$, then $F_{a_1}(x)\geq F_{a_2}(x)$
for every $x\in [-1,1]$. In particular, for every $a>0$,
$$
F_a(x)\geq F_0(x)=T_n(x)+x+1-\frac{2x+1}{n^2}\,T_n^{\prime}(x)
=\varphi_n(x)\,,\qquad x\in [-1,1]\,.
$$
In view of Theorem~\ref{t2}, $\varphi_n(x)\geq 0$ for every $x\in
[0,1]$, with the equality holding only for $x=1$ and if $n=4k+2$,
for $x=0$. Therefore, with $a$ as given in \eqref{e4.1}, we have
$F_a(x)\geq \varphi_n(x)>0$ for every $x\in (0,1)$, and it remains
to prove inequality \eqref{e4.2} and clarify the cases of equality
only when $x\in [-1,0]$.

Let us denote by
$$
x_k=\cos\frac{k\pi}{n},\qquad k=0,\ldots,n\,,
$$
the zeros of $(1-x^2)T_n^{\prime}(x)$, and let $H(f;x)$ be the
Hermite interpolating polynomial, with interpolation nodes
$x_0,x_1,x_1,x_2,x_2,\ldots,x_{n-1},x_{n-1}, x_n$, for a
differentiable function $f$, i.e., $H(f;x)$ is determined by the
conditions
\begin{eqnarray*}
&&H(f;x_k)=f(x_k),\quad k=0,1,\ldots,n,\\
&&H^{\prime}(f;x_k)=f^{\prime}(x_k), \quad k=1,\ldots,n-1\,.
\end{eqnarray*}
A straightforward calculation shows that
\begin{equation}\label{e4.3}
\begin{split}
H(f;x)=&\frac{\big[T_n^{\prime}(x)\big]^2}{2n^4}\,\big[(1+x)\,f(x_0)
+(1-x)\,f(x_n)\big]\\
&+(1-x^2)\,\sum_{k=1}^{n-1} \frac{\ell_k^2(x)}{(1-x_k^2)^2}\,
\mathcal{L}_k(f;x)\,,
\end{split}
\end{equation}
where, for $k=1,\ldots, n-1$, $\ell_k$ are the Lagrange basis
polynomials for interpolation at the zeros of $T_n^{\prime}$,
$$
\ell_k(x)=\frac{T_n^{\prime}(x)}{(x-x_k)T_n^{\prime\prime}(x_k)},
$$
and
\begin{equation}\label{e4.4}
\mathcal{L}_k(f;x):=(1-x_k\,x)f(x_k)
+(1-x_k^2)(x-x_k)f^{\prime}(x_k)\,.
\end{equation}

It follows from the uniqueness of the Hermite interpolation
polynomial that $H(f;\cdot)\equiv f(\cdot)$ whenever $f$ is a
polynomial of degree at most $2n-1$, in particular,
$H(f_i;\cdot)\equiv f_i(\cdot)$ for $i=1,\,3$ and
$H(F_a;\cdot)\equiv F_a(\cdot)$\,.

From $T_n(x_k)=(-1)^k$, $0\leq k\leq n$, and $T_n^{\prime}(x_k)=0$,
$1\leq k\leq n-1$, $T_n^{\prime}(x_0)=n^2$,
$T_n^{\prime}(x_n)=(-1)^n\,n^2$, we find
\begin{equation}\label{e4.5}
\begin{split}
& f_1(x_0)=f_3(x_0)=0,\\
&f_1(x_k)=2+(-1)^k,\ \ f_3(x_k)=1-x_k,\quad 1\leq k\leq n-1;\\
&f_1(x_n)=f_3(x_n)=2(1+(-1)^n)\,.
\end{split}
\end{equation}
In particular, \eqref{e4.5} yields
$$
F_a(x_0)=0\,,\quad F_a(x_n)=2a\big[1+(-1)^n\big]\,,
$$
in agreement with the claim of Theorem~\ref{t3} that $F_a(x)$
vanishes at $x_0$ and if $n$ is odd, at $x_{n}$; in addition,
$F_a(x_n)>0$ if $n$ is even and $a>0$. Moreover, from
$F_a(\cdot)\equiv H(F_a;\cdot)$, \eqref{e4.3} and \eqref{e4.4} we
infer
\begin{equation}\label{e4.6}
F_a(x)=\frac{a\big[1\!+\!(\!-\!1)^{n}\big]}{n^4}\,
(1\!-\!x)\big[T_n^{\prime}(x)\big]^2
\!+\!(1\!-\!x^2)\,\sum_{k=1}^{n-1}\frac{\ell_k^2(x)}{(1\!-\!x_k^2)^2}
\,\mathcal{L}_k(F_a;x)\,.
\end{equation}
In order to find $\mathcal{L}_k(F_a;x)$, $1\leq k\leq n-1$\,, we
firstly evaluate $\mathcal{L}_k(f_i;x)$, $i=1,\,3$\,.

Making use of the differential equation \eqref{e2.2} and
$$
f_1^{\prime}(x)=\Big(1-\frac{1}{n^2}\Big)\,T_n^{\prime}(x)
-\frac{x+2}{n^2}\,T_n^{\prime\prime}(x), \quad
f_3^{\prime}(x)=\frac{1}{n^2}\,
\big[T_n^{\prime}(x)-(1-x)\,T_n^{\prime\prime}(x)\big]-1\,,
$$
we find
\begin{equation}\label{e4.7}
f_1^{\prime}(x_k)=(-1)^{k}\,\frac{x_k+2}{1-x_k^2}\,,\quad
f_3^{\prime}(x_k)=\frac{(-1)^{k}}{1+x_k}-1\,,\qquad
k=1,\ldots,n-1\,.
\end{equation}
From \eqref{e4.4}, \eqref{e4.5} and \eqref{e4.7} we obtain
$$
\mathcal{L}_k(f_1;x)=(1-x_k\,x)\big(2+(-1)^k\big)
+(-1)^{k}(x_k+2)(x-x_k)\,,
$$
or, equivalently,
\begin{equation}\label{e4.8}
\mathcal{L}_k(f_1;x)=
\begin{cases}(1-x_k)(3+2x+x_k),& k \text{ - even},
\vspace*{1ex}\\
(1+x_k)(1+x_k-2x),& k \text{ - odd}.
\end{cases}
\end{equation}
Likewise, we get
$$
\mathcal{L}_k(f_3;x)=(1-x_k\,x)(1-x_k)
+(x-x_k)\,\big[(-1)^{k}(1-x_k)-1+x_k^2\big]\,,
$$
which simplifies to
\begin{equation}\label{e4.9}
\mathcal{L}_k(f_3;x)=
\begin{cases}(1-x_k)(1+x_k^2-2x_k\,x),& k \text{ - even},
\vspace*{1ex}\\
(1-x_k^2)(1+x_k-2x),& k \text{ - odd}.
\end{cases}
\end{equation}
Now $\mathcal{L}_k(F_a;x)=(1+a)\,\mathcal{L}_k(f_1;x)
-\mathcal{L}_k(f_3;x)$, \eqref{e4.8} and \eqref{e4.9} yield
\begin{proposition}\label{pp2.5}
For $1\leq k\leq n-1$, we have
$$
\mathcal{L}_k(F_a;x)=
\begin{cases}
(1-x_k)\big[2(1+x)(1+a+x_k)+(a-x_k)(1+x_k)\big],& k \text{ - even},
\vspace*{1ex}\\
(a+x_k)(1+x_k)(1+x_k-2x),& k \text{ - odd}.
\end{cases}
$$
\end{proposition}

Let us note that \eqref{e4.6} and Proposition~\ref{pp2.5} hold for
an arbitrary constant $a$. We show below that, with $a$ as given in
\eqref{e4.1}, $\mathcal{L}_k(F_a;x)>0$ for $x\in (-1,0]$, except for
a specific index $k$. When $k$ is even, we prove even more.
\begin{proposition}\label{pp3}
If $n\ge 4$ is even, $k$ is even, $1<k< n-1$, and $a=x_1$, then
$$
\mathcal{L}_k(F_a;x)>0,\qquad -1\leq x\leq 1\,.
$$
\end{proposition}

\proof In this case $a-x_k=x_1-x_k>0$, and it follows from
Proposition~\ref{pp2.5} that
$$
\mathcal{L}_k(F_a;x)> 2(1+x)(1-x_k)(1+a+x_k)\geq 0,\qquad -1\leq
x\leq 1\,.
$$
\qed

\begin{proposition}\label{pp4}
If $n\geq 4$, $k$ is odd, $1\leq k\leq n-1$, and $a$ is given by
\eqref{e4.1}, then
\begin{equation*}
\mathcal{L}_k(F_a;x)\begin{cases} >0, & -1<x\leq 0,\quad 1\leq
k<n-2\,, \vspace*{1ex}\\
\equiv 0,& k=n-2,\ \ n\ \text{ - odd},\vspace*{1ex}\\
\equiv 0,& k=n-1,\ \ n\ \text{ - even}\,.
\end{cases}
\end{equation*}
\end{proposition}
\proof If $n$ is even, then $a=x_1=-x_{n-1}$, thus $a+x_{n-1}=0$ and
$\mathcal{L}_{n-1}(F_a;x)\equiv 0$ in view of
Proposition~\ref{pp2.5}. If $n$ is odd, then $a=x_2=-x_{n-2}$ and
$a+x_{n-2}=0$, which by Proposition~\ref{pp2.5} implies
$\mathcal{L}_{n-2}(F_a;x)\equiv 0$. Finally, if $k<n-2$ is odd, then
$a+x_k>a+x_{n-2}\geq x_2+x_{n-2}=0$, and from
Proposition~\ref{pp2.5} we infer
$$
\mathcal{L}_k(F_a;x)=(a+x_k)(1+x_k)(1+x_k-2x)>0, \qquad -1\leq x<
\frac{1+x_{n-3}}{2}.
$$
Since $(1+x_{n-3})/2>0$, it follows that $\mathcal{L}_k(F_a;x)>0$
for every $x\in [-1,0]$. \qed
\medskip

We are ready to accomplish the proof of the claim of
Theorem~\ref{t3} in the case $x\in [-1,0]$.\smallskip

Assume first that $n\geq 4$ is even. If $a=x_1$, then \eqref{e4.6}
and Proposition~\ref{pp4} imply
\begin{equation}\label{e4.10}
F_{x_1}(x)=\frac{2x_1}{n^4}\,(1-x)\big[T_n^{\prime}(x)\big]^2 +
(1-x^2)\,\sum_{k=1}^{n-2}\frac{\ell_k^2(x)}{(1-x_k^2)^2}
\,\mathcal{L}_k(F_a;x)\,.
\end{equation}
Since $\ell_k(x_{n-1})=0$ for $1\leq k\leq n-2$ and
$T_n^{\prime}(x_{n-1})=0$, it follows that
$F_{x_1}(x_{n-1})=F_{x_1}(-x_1)=0$. If, on the other hand, $x\in
(-1,0]$ and $x\ne x_{n-1}$, then all the summands in the sum in the
right-hand side of \eqref{e4.10} are non-negative, by virtue of
Propositions~\ref{pp3} and \ref{pp4}, with at least one of them
strictly positive, therefore $F_{x_1}(x)>0$ in this case. If
$a<x_1$, then from \eqref{e4.6} and Proposition~\ref{pp2.5} we
obtain
\begin{equation*}
\begin{split}
F_a(x_{n-1})&=(1-x_{n-1}^2)\,\sum_{k=1}^{n-1}
\frac{\ell_k^2(x_{n-1})}{(1-x_k^2)^2} \,\mathcal{L}_k(F_a;x_{n-1})\\
&=\frac{\mathcal{L}_{n-1}(F_a;x_{n-1})}{1-x_{n-1}^2}
=a+x_{n-1}<x_1+x_{n-1}=0\,,
\end{split}
\end{equation*}
showing that the inequality $F_a(x)\geq 0$ fails to hold for
$x=x_{n-1}$.\smallskip

Now assume that $n\geq 5$ is odd. If $a=x_2$, then \eqref{e4.6} and
Proposition~\ref{pp4} imply
\begin{equation}\label{e4.11}
F_{x_2}(x)=
(1\!-\!x^2)\,\Big[\sum_{k=1}^{n-3}\frac{\ell_k^2(x)}{(1-x_k^2)^2}
\,\mathcal{L}_k(F_a;x)\!+\!\frac{\ell_{n-1}^2(x)}{(1-x_{n-1}^2)^2}
\,\mathcal{L}_{n-1}(F_a;x)\Big]\,.
\end{equation}
Since $\ell_k(x_{n-2})=0$ for $k\ne n-2$, we have
$F_{x_2}(x_{n-2})=0$. If $x\in (-1,0]$, $x\ne x_{n-2}$, then, in
view of Propositions~\ref{pp3} and \ref{pp4}, all the summands in
the brackets in \eqref{e4.11} are non-negative, and at least one of
them is strictly positive, therefore $F_{x_2}(x)>0$. Finally, if
$a<x_2$, then from \eqref{e4.6} and Proposition~\ref{pp2.5} we find
\begin{equation*}
\begin{split}
F_a(x_{n-2})&=(1-x_{n-2}^2)\,\sum_{k=1}^{n-1}
\frac{\ell_k^2(x_{n-2})}{(1-x_k^2)^2} \,\mathcal{L}_k(F_a;x_{n-1})\\
&=\frac{\mathcal{L}_{n-2}(F_a;x_{n-2})}{1-x_{n-2}^2}
=a+x_{n-2}<x_2+x_{n-2}=0\,,
\end{split}
\end{equation*}
i.e., $F_a(x_{n-2})<0$ if $a<x_2$.

The proof of Theorem~\ref{t1} is complete.

\noindent
{\sc Geno Nikolov} \smallskip\\
Department of Mathematics and Informatics\\
Universlty of Sofia \\
5 James Bourchier Blvd. \\
1164 Sofia \\
BULGARIA \\
{\it E-mail:} {\tt geno@fmi.uni-sofia.bg}

\end{document}